\documentclass{article}
\usepackage{amsfonts}
\usepackage{amsmath}
\usepackage{graphicx}

\newtheorem{theorem}{Theorem}[section]
\newtheorem{corollary}{Corollary}[section]
\newtheorem{definition}{Definition}[section]

\newtheorem{proposition}[theorem]{Proposition}

\begin{document}

\title{Can a fingerprint be modelled by a differential equation ?}
\author{F. Zinoun \\
LabMIA-SI, Department of Mathematics, Faculty of Sciences, University
Mohammed V, Rabat, Morocco }

\begin{center}
{\Huge Can a Fingerprint be Modelled by a Differential Equation ?}

\bigskip

\bigskip

\textit{(When Galton Meets Poincar\'{e})}

\bigskip

\bigskip

\bigskip

by

Fouad ZINOUN

\bigskip

\textit{\small LabMIA, Department of Mathematics, Faculty of Sciences}

\textit{\small Mohammed V University in Rabat - Morocco}

\textit{\small zinoun@fsr.ac.ma}
\bigskip
\bigskip
\bigskip
\bigskip

\includegraphics[scale=0.7]{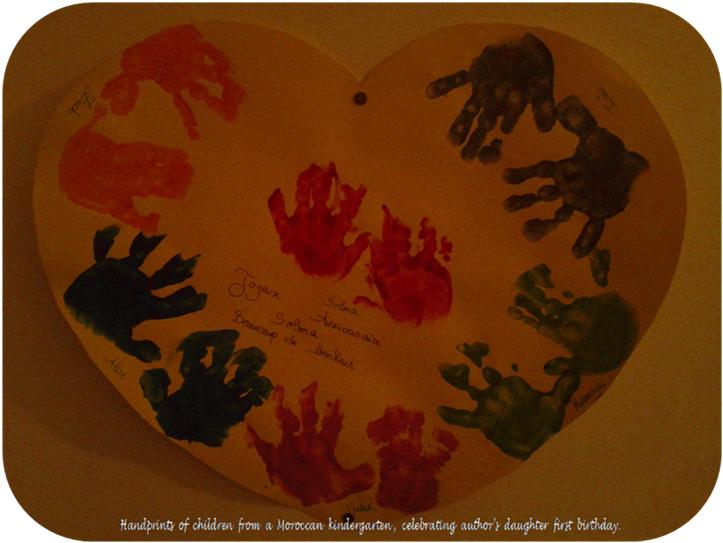}


\end{center}

\tableofcontents
\textbf{References}  
{\it (Let's Twist Again, a parody song)}

\newpage

\textbf{Abstract.}

Some new directions to lay a rigorous mathematical foundation for the phase-portrait-based modelling of fingerprints are discussed in the present work. Couched in the language of dynamical systems, and preparing to a
preliminary modelling, a back-to-basics analogy between Poincar\'{e}'s
categories of equilibria of planar differential systems and the basic
fingerprint singularities according to Purkyn\v{e}-Galton's standards is
first investigated. Then, the problem of the global representation of a
fingerprint's flow-like pattern as a smooth deformation of the phase
portrait of a differential system is addressed. Unlike visualisation in
fluid dynamics, where similarity between integral curves of smooth vector
fields and flow streamline patterns is eye-catching, the case of an oriented
texture like a fingerprint's stream of ridges proved to be a hard problem
since, on the one hand, not all fingerprint singularities and nearby
orientational behaviour can be modelled by canonical phase portraits on the
plane, and on the other hand, even if it were the case, this should lead to
a perplexing geometrical problem of connecting local phase portraits, a
question which will be formulated within Poincar\'{e}'s index theory and
addressed via a normal form approach as a bivariate Hermite interpolation
problem. To a certain extent, the material presented herein is
self-contained and provides a baseline for future work where, starting from
a normal form as a source image, a transport via large deformation flows is
envisaged to match the fingerprint as a target image.

\section{Introduction}

In 1892, Sir Francis Galton, the English Victorian scientist, published the
first book of his \textquotedblleft trilogy\textquotedblright\ on
fingerprints \cite{Galton}. As said by the author, his attention had been
first drawn to ridges when preparing, some years ago, a Royal Institution
lecture on personal identification, which aimed at an account of the newly
introduced French anthropometric method of Alphonse Bertillon \cite{Bert}.
Realising both how much had been done on the subject and how much there
remained to do, and being chiefly based on a thesis of the Czech
physiologist, Jan Purkyn\v{e}, at the university of Breslau \cite{Purk} - a
very rare pamphlet on classification of papillary ridges, he will become
perhaps the first to place the fingerprint-based William Herschel's
identification method \cite{Hers} on a scientific footing and to lay
securely the foundation of a new branch of inquiry. Following a series of
memoirs upon the subject [5-8], a system for classifying fingerprint patterns into three broad
categories, which is very useful for rough preliminary purposes, is then
mainly used in Galton's book and of which frequent reference will be made in
this paper: Arches, Loops and Whorls (ALW), in a sense to be specified.
\bigskip

In the same year, Henri Poincar\'{e}, the French universal scientist and the
father of the qualitative theory of differential equations, published the
first volume of his \textquotedblleft trilogy\textquotedblright\ on
celestial mechanics \cite{Poincare}, a masterpiece written during the last
decade of the nineteenth century, following the pioneering works of Cauchy,
Lagrange and Laplace, his own inaugural thesis \cite{Poincare1} and a series
of papers \cite{Poincare2} where different types of singularities have been
named and studied: Nodes, Foci, Saddles and Centers.

\bigskip
\bigskip

\begin{figure}
\begin{center}
\includegraphics{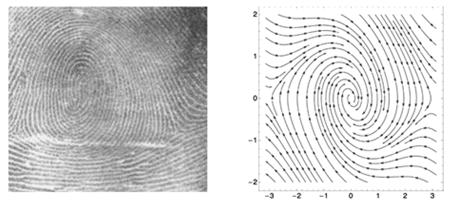} 
\end{center}
\caption{\it A whorl fingerprint
pattern (NIST Database) and the computer-drawn phase
portrait of a simple damped pendulum.}
\label{basic}
\end{figure}

The two scientists have probably never met in person,\footnote{%
If one has to establish a far-fetched link\ between Galton and Poincar\'{e},
it would be perhaps a Poincar\'{e}'s work of 1885 on the equilibrium figures
of a fluid mass \cite{Poincare3}, from which George Darwin, the son of
Charles, who was Galton's half-cousin, deduced what he believed to be a
mechanism for the formation of the Moon!} but like Monsieur Jourdain who
was speaking prose without knowing it, they were perhaps speaking the same
language, as insinuated in Fig.\ref{basic}: fingerprinting when dressing
phase portraits for Poincar\'{e}, and conversely, for Galton, solving
differential equations when deciphering fingerprints!

\bigskip

The idea is then the following: to what extent can a fingerprint's
orientation image be visualised as (a smooth deformation of) the phase
portrait (or more faithfully, in Poincar\'{e}'s language, a system of
characteristics) of a planar dynamical system? In other terms, to what
extent can a classification system of fingerprints be couched in the
language of the qualitative theory of differential equations?

\bigskip

Contrary to what this introduction may be suggesting, the idea of using
phase portraits in texture modelling is not new, as can be already seen
through an interdisciplinary program initiated by the Semiconductor Research
Corporation at the University of Michigan in the 80's, whose goal was to
develop a visual language for representing visual data in semiconductor
wafer processing. The thesis by Rao \cite{Rao}, for instance, was part of a
larger effort within the program to device a symbolic description of
oriented texture patterns using the qualitative theory of differential
equations. Curiously, a large album of real texture images have been
analysed, like an invigorating wood knots accord with notes of orange peel,
enhanced by hints of brush strokes from Van Gogh paintings, but no
fingerprint image seems to have been considered. Special mention must also
be made of the thesis by Ford \cite{Ford}, at the University of Arizona, on
2-D fluid flow modelling and visualisation, whereby complex flows were split
into simpler and easily described components, the latter being modelled by
linear phase portraits and then combined to obtain a global model for the
entire flow field. This idea has been finally applied to fingerprints by Li
and Yau and recapitulated as a chapter in \cite{LY} (see also references
therein) where, following Kass and Witkin's scheme of squaring the gradient
vectors in the computed oriented texture field \cite{KW}, basic fingerprint
singular points have been assimilated to either a Focus or a Saddle.
However, amid the wealth of literature on fingerprints modelling in the last
two decades, if we restrict ourselves solely to the Taylor expansion as a
basis for phase portrait models,\footnote{%
For other variations, see for instance Wang et al. \cite{WHP} for a
Fourier-expansion-based approach and Ram et al. \cite{RBB} for a
Legendre-polynomial-based phase portrait model.} and while this clearly
yields good results at a practical level, there are still some outstanding
issues to be mathematically clarified. In fact, at a purely mathematical
level, and without going into details for the moment, previous
phase-portrait-based methods will be most certainly found to work near a
singular point for a fingerprint orientation as long as the Hartman-Grobman
theorem works for a nonlinear dynamical system near an equilibrium point.
This can also be expressed in terms of Peixoto's theorem within the scope of
structural stability. Indeed, when dealing for example with a Saddle or a
Focus, we know that a polynomial perturbation of the linear terms will not
change the nature of these points (see the discussion in paragraph 5.2.1 on
robustness of singular points against addition of nonlinear terms). However,
for a Center, which is an unavoidable point in the Circular/Elliptical Whorl
or the Circlet in Loop modelling (see below), let alone a degenerate point
(see the classification in paragraph 3.2), it is mathematically possible,
but highly unlikely, that a center-like behaviour near a singular point will
be preserved after addition of nonlinear terms, unless appropriate
symmetries are shown by the global flow, a behaviour which is often observed
in natural occurring flows, but certainly not on fingerprints. Besides, only
nondegenerate singular points have been considered for the derived piecewise
phase portrait models, then going to miss (nontrivial) degenerate points
occurring in interesting bifurcations, next to the rich variety of phase
portraits they offer (see Fig.\ref{singular} below). Just in the
nondegenerate case, a large part has been in fact devoted to the theory of
centers by Poincar\'{e} [11, chap.XI], to show that the
question has to be considered in its own right. At another level, the
problem of interconnecting singular points should not be dealt with as if it
were always solvable; obviously, this is true in fluid dynamics where
connecting streamlines can be experimentally visualised and then
topologically modelled. However, for a texture image like a fingerprint, and
unless not more than a restricted and well-chosen number of relevant
singular points are considered, such a connexion is generally impossible as
it comes out from a global index theory analysis of smooth vector fields on
two-dimensional surfaces.

\bigskip

The present paper has no pretension to outperform previous works on
phase-portrait-based modelling, but a targeted approach will be adopted to
mobilising sophisticated results from dynamical systems, in line with Poincar%
\'{e}'s \oe uvre, to help overcome the limited capability of a basic
linearisation approach and open up new possibilities for future works on
fingerprints modelling. Nonhyperbolic and degenerate cases will be then
carefully handled within a structure-preserving normal form approach, taking
profitably the large variety of qualitative behaviours they host. As will be
shown, if the difficulty can be partially overcome for a preliminary
modelling within the ALW classification, in which no more than three
singular points are taken into account, the general case involving many
singularities with different natures, however, leads to a serious
mathematical problem. Explicitly, even if all basic fingerprint patterns
were modelled by appropriate phase portraits, the problem of associating a
global phase portrait to a fingerprint would lead to an advanced bivariate
Hermite interpolation problem for which an algebraic solution in the general
case is hopeless. To obtain some elements of answer, the geometrical notion
of connecting local phase portraits will be considered from an intuitive
point of view, in the spirit of Poincar\'{e}'s question of how singular
points of vector fields are distributed in the phase space, and how the
study of a function defined in the vicinity of a singular point can be
extended to the whole space. The geometrical notion of connecting local
phase portraits, however, should be much easier to address than a direct
attack of the interpolation problem. It will be then easily seen, within
Poincar\'{e}'s index theory, that such a connexion do not always exist, and
when it exists, obviously, it is not necessarily unique. But in general,
there is no systematic approach to carry out a connexion, although a lot has
been done by Poincar\'{e} on the sphere, after gnomonic projection, from
which the subject can be shown to derive some strength and fruitfulness.

\bigskip

Globally, the paper is structured as follows. In the first part, a brief survey of Galton's book is first given, followed by a
(partial, yet sufficient) classification of singular points on the plane. A
gallery of local phase portraits is then presented from which a collection
of typical singularities will be hand-picked to match - whenever possible -
the basic fingerprint patterns to be modelled. Then, an attempt at a
preliminary modelling of fingerprints is made according to Purkyn\v{e}%
-Galton's standards, where some hand-drawn phase portraits are given as
approximations. As will be seen, the main obstruction caused by a delta-like
pattern will be overcome by integrating in the phase plane a cusp-like
singular point from the well-known Bogdanov-Takens bifurcation. In the
second part, some directions for building global phase portraits
from local ones are discussed within Poincar\'{e}'s index theory and, when
possible, a structure-preserving normal form\ - in a sense to be specified -
will be computed as a bivariate Hermite interpolation problem and attributed
to some subclasses within the ALW system. As said in the abstract, such a
normal form will be the starting point for future work where, to obtain a final
signature, a transport via large deformation flows is envisaged to carry
away the normal form as a source image in order to match the fingerprint as
a target image.

\part{\protect\Large Preliminary phase-portrait-based modelling of
fingerprints}

\section{A brief survey of Galton's book}

\begin{figure}
\begin{center}
\includegraphics{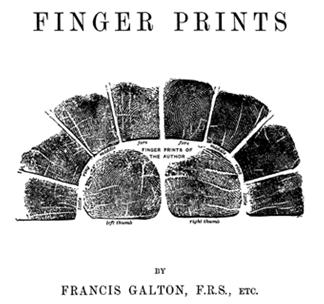} 
\end{center}
\caption{\it Title-page of Galton's book with author's fingerprints.}
\label{}
\end{figure}


In [1], and not without a wicked sense of humour, Galton begins with these
words: \textit{\textquotedblleft The palms of the hands and the soles of
the feet are covered with two totally distinct classes of marks. The most
conspicuous are the creases or folds of the skin which interest the
followers of palmistry, but which are no more significant to others than the
creases in old clothes; they show the lines of most frequent flexure, and
nothing more. The least conspicuous marks, but the most numerous by far, are
the so-called papillary ridges; they form the subject of the present book.
If they had been only twice as large as they are, they would have attracted
general attention and been commented on from the earliest times. Had Dean
Swift known and thought of them, when writing about the Brobdingnags, whom
he constructs on a scale twelve times as great as our own, he would
certainly have made Gulliver express horror at the ribbed fingers of the
giants who handled him. The ridges on their palms would have been as broad
as the thongs of our coach-whips.\textquotedblright }

\bigskip

\begin{figure}
\begin{center}
\includegraphics[scale=0.8]{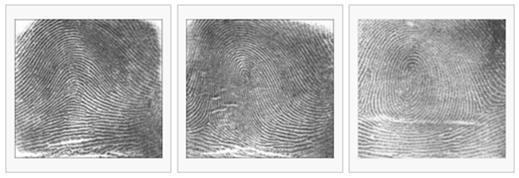} 
\end{center}
\caption{\it From left to right: Arch, Loop and Whorl (NIST Database).}
\label{ALW}
\end{figure}

After treating of the previous use of fingerprints, from superstition of
personal contact to the modern regular official employment of Herschel,
various methods of making good prints and enlarging them are described as
they are adopted at the author's own anthropometric laboratory. Then these
preliminary topics having been disposed of, the author begins with a
discussion of the various patterns formed by the lineations, emphasing the
independent ones that appear upon the bulbs of the fingers, and where plates
of the principal varieties of patterns are given as a visual support. A useful
classification system for rough preliminary study is then presented into
Arches, \textit{for which we have no pattern strictly speaking}; Loops, 
\textit{where we have a system of ridges that bends back upon itself and in
which no one ridge turns through a complete circle}; Whorls, \textit{for
which at least one ridge turns through a complete circle }(Fig.\ref{ALW}).
Of course, chapters dealing with evidential values, methods of indexing,
personal identification or heredity are beyond the interest of the present
paper.

\bigskip

An important passage of the book should be however highlighted: a
translation in part from the Latin of the Commentatio of Purkyn\v{e}, made
by the author himself, which a copy has been procured from the United States
to the Library of the Royal College of Surgeons. The following nomenclature
was then established according to the nine principal varieties of curvature
observed by Purkyn\v{e}, and presented here in the same order as they appear
in Plate 12 of Galton's book (see Fig.\ref{Purk}): 1) Transverse flexures;
2) Central Longitudinal Stria; 3) Oblique Stria; 4) Oblique Sinus; 5)
Almond; 6) Spiral; 7) Ellipse, or\ Elliptical Whorl; 8) Circle, or\ Circular
Whorl; 9) Double Whorl. Following Purkyn\v{e}, all these forms have been
concisely described by Galton, within Arches for 1-3, Loops for 4-5 and
Whorls for 6-9. We prefer to refer to diagrams for explanation at this
stage, while more detail will be given for each variety when proceeding to
modelling.

\begin{figure}
\begin{center}
\includegraphics{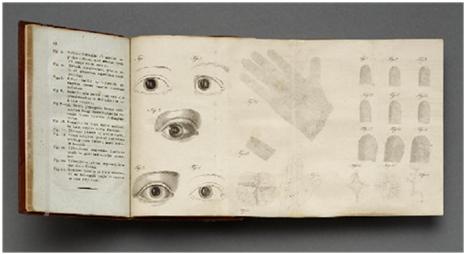} 
\end{center}
\caption{\it Purkyn\v{e}'s Commentatio (U.S. National Library of Medicine).}
\label{}
\end{figure}


\begin{figure}
\begin{center}
\includegraphics{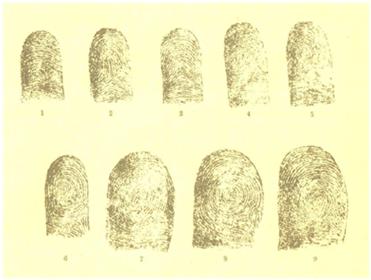} 
\end{center}
\caption{\it The standard patterns of Purkyn\v{e}, as captured from Galton's book.}
\label{Purk}
\end{figure}


\bigskip

To start the study in section 4, and for purely aesthetic reasons, I'll
choose from the above configurations the Spiral and the Elliptical/Circular
Whorl as basic models to be identified to the phase portrait of an
autonomous planar system of ordinary differential equations. For such a
purpose, I'll try to establish an analogy between some categories of
singular points of a planar vector field and the basic singular points\
composing the fingerprint's orientation image being analysed. Reference is
made to Galton's book (and mainly to the translation of Purkyn\v{e}'s thesis
therein) when it comes to using terminology from fingerprints, except for
the term \textit{Delta} which is borrowed from the well-known Henry
Classification System \cite{Henry}\footnote{%
Expanding on Galton's classification system, this is another interesting
book which, by order of the government of (British) India, and for
bureaucratic settings, was published in 1900 by Edward Henry as being a
former member of the civil service at the presidency of Fort William in
Bengal.} and will be preferred to Purkyn\v{e}-Galton's \textit{Triangle}. As
for terminology from dynamical systems, reference will be mainly made to
Poincar\'{e}, whose the \oe uvre is the basis of the planar classification
of singular points given below.

\section{Poincar\'{e}'s imprint}

In July 2012, celebrating 100 years after Poincar\'{e}, Alain Chenciner,
from the Institut de M\'{e}canique C\'{e}leste de l'Observatoire de Paris,
delivered in the cemetery of Montparnasse a moving and eloquent speech, with
a riveting and poetic survey of all Poincar\'{e}'s \oe uvre. Besides the man
and his heritage, one of the interesting things to which attention could be
drawn in Chenciner's speech, is especially an \textquotedblleft
opposition\textquotedblright\ of Poincar\'{e}'s spirit and nature to Charles
Hermite's, the French realistic and anti-geometer mathematician. In fact,
within the present paper, the way in which is addressed the problem of
assigning (a smooth deformation of) a differential equation's phase portrait
to a fingerprint, i.e. whether via the intuitive notion of connecting local
phase portraits or, rigorously, as an algebraic multivariate interpolation
problem, could be seen as a certain expression of the opposition between
these two natures, and will be sometimes perceptible as work progresses.

\subsection{Singular points of the first species}

Following [9-11], let consider curves defined by
an equation of the form%
\begin{equation}
\frac{dx}{X}=\frac{dy}{Y}  \label{1}
\end{equation}%
where $X$ and $Y$ are analytic functions in $x$ and $y.$ Such curves are
called \textit{characteristics} by Poincar\'{e}. As we are not concerned
with the study of infinite branches for the moment, we will not consider the
gnomonic projection on the sphere as Poincar\'{e} usually does, restricting
ourselves to (a subdomain of) the phase plane $(x,y).$ Following Cauchy \cite%
{Cauchy}, Briot and Bouquet \cite{BB}, and himself \cite{Poincare1}, Poincar%
\'{e} gave a complete description of the characteristics near an isolated
nondegenerate singular point, that is for which $X=Y=0$ and the Jacobian
matrix $J$ (of first-order partial derivatives of $X$ and $Y$ at the
singular point) has no zero eigenvalues. Depending on the distribution of
these eigenvalues in the complex plane, a classification into four
categories is given: Nodes, Foci, Saddles, for the hyperbolic case (no
eigenvalues with zero real part), and Centers, for the nonhyperbolic case.
Poincar\'{e} calls them singular points of the first species. In Fig.\ref%
{diagram}, a diagram of bifurcation is given according to the trace $\tau $
and the determinant $\triangle $ of the Jacobian matrix, the half-axis $\tau
=0,\triangle >0$ and the axis $\triangle =0$ corresponding to nonhyperbolic
singularities that arise at Andronov-Hopf and Saddle-Node Bifurcation,
respectively.

\begin{figure}
\begin{center}
\includegraphics[scale=0.85]{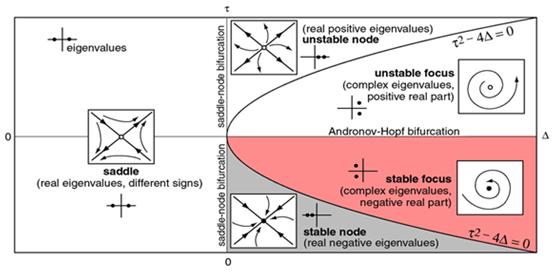} 
\end{center}
\caption{\it Classification of
singular points of planar systems (E. Izhikevich, Scholarpedia (2007)).}
\label{diagram}
\end{figure}


\subsection{Singular points of the second species}

The case when the curves $X=Y=0$ intersect in many combined points leads to
singular points of the second species, which are indeed nonhyperbolic and
can be considered as the limit of a system of singular points of the first
species, combined together. They are sometimes referred to as multiple
singular points since they can be made to split into a number of hyperbolic
critical points under suitable perturbation of $X$ and/or $Y.$

\bigskip

\begin{figure}
\begin{center}
\includegraphics[scale=0.8]{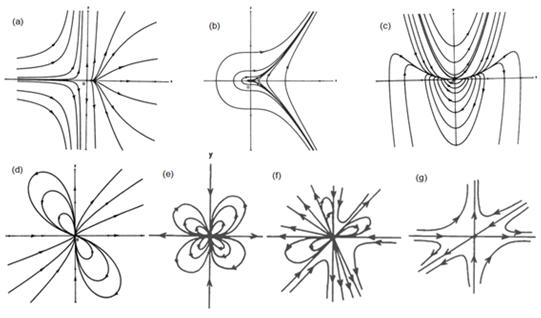} 
\end{center}
\caption{\it Some degenerate
singular points, as collected from Perko's book.$^4$}
\label{singular}
\end{figure}

As said by Poincar\'{e} \cite{Poincare2}\ in p. 393, such points are of too
numerous and too diverse particularities to be studied in details. In Perko 
\cite{Perko}, following Poincar\'{e} \cite{Poincare2}, Bendixson \cite%
{Bendixson} and more recently Andronov et al. \cite{Andronov}, a collection
of interesting results on nonhyperbolic singular points of planar analytic
systems is recalled and a gallery of phase portraits is plotted for
different singular points with sometimes unusual behaviour of nearby
trajectories (Fig.\ref{singular}).

\bigskip

\begin{figure}
\begin{center}
\includegraphics[scale=0.9]{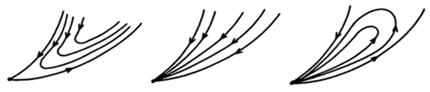} 
\end{center}
\caption{\it From left to right, a
hyperbolic, a parabolic and an elliptic sector \protect\cite{Perko}.}
\label{sector}
\end{figure}

In the case when $J$ has at least one zero eigenvalue (the degenerate case),
but $J\neq 0$, it is shown that there are at most $2(m+1)$ directions $%
\theta $ along which a solution curve of $(1)$ may approach the singular
point (supposed put at the origin, without loss of generality), provided the
function
{\let\thefootnote\relax\footnote{$^4$May I ever be excused for \textquotedblleft stealing\textquotedblright\ these phase portraits from \cite{Perko}, without permission from Prof. Lawrence Perko.}}%

\begin{equation*}
f(\theta )=\cos \theta X_{m}(\cos \theta ,\sin \theta )-\sin \theta
Y_{m}(\cos \theta ,\sin \theta )
\end{equation*}%
is not identically zero, $X_{m}$ and $Y_{m}$ being the $m$th-degree terms
with which begin the Taylor series of $X$ and $Y,$ respectively. These
directions are given by solutions of the equation $f(\theta )=0$ and then
the notion of \textit{sector} become fundamental for the classification. In
fact, a sufficiently small neighborhood of the origin will be divided by
these curves into a finite number of open regions (sectors), each of them
being either of a hyperbolic, a parabolic or an elliptic type (Fig.\ref{sector}). This is to be
understood in a topological sense, no regard being paid to the direction of
the flow.

\bigskip

The trajectories which lie on the boundary of a hyperbolic sector are then
called separatrices, and it is shown that, besides the usual types of
singular points for planar analytic systems, i.e. Nodes, Foci, Saddles and
Centers, the only other types of singular points that can occur for $(1)$
(when $J\neq 0$) are: Saddle-Nodes, for which we have two hyperbolic sectors
and one parabolic (Fig.\ref{singular}-a); Cusps, occurring in
Bogdanov-Takens bifurcation and for which we have two (and only two)
hyperbolic sectors (Fig.\ref{singular}-b); Singular points with an elliptic
domain, for which we have one elliptic sector, one hyperbolic sector and two
parabolic sectors (Fig.\ref{singular}-c).

\bigskip

More precisely, if $J$ has exactly one zero eigenvalue, the underlying
singular point is either a Node, a Saddle, or a Saddle-Node; and if $J$ has
two zero eigenvalues (remember, $J\neq 0),$ the singular point is either a
Focus, a Center, a Node, a Saddle, a Saddle-Node, a Cusp, or a singular
point with an elliptic domain. Note that we have only taken into account
analytic vector fields, so according to a theorem by Dulac (1923) on the
finitude of the number of limit cycles (i.e. isolated closed trajectories, a
question intimately linked to the famous and still unsolved Hilbert's
sixteenth problem), the case of the Center-Focus is to be excluded: this
corresponds to a singular point surrounded by an infinite number of
accumulating limit cycles. As will be seen, however, this is not a big loss
since patterns of a spiral or a circular form can be merely represented by a
Focus or a Center.

\bigskip

The case $J=0$ is a bit exceptional, and the behaviour near the singular
point can be more complex (Fig.\ref{singular}-(d-g)). All we recall here is
a (consequence of a) beautiful result due to Bendixson \cite{Bendixson},
within index theory (cf. subsection 5.1.1), and which can be summarised as
follows:%
\begin{equation*}
e=h \pmod 2
\end{equation*}%
$e$ and $h$ being, respectively, the number of elliptic sectors and the
number of hyperbolic ones. In other terms, $e$ and $h$ have the same parity.
Note that all degenerate singular points (a-g) in Fig.\ref{singular} agree
well with this formula. One can of course endeavour to enlarge our
collection by including other exotic points in the degenerate case; this
would not be however of much use since, as we shall see in the next section,
the formula above is quite sufficient to test whether an imaginary singular
point (suggested by a pattern) can match a real singular point or not.

\section{Preliminary modelling}

In this section, we will consider Purkyn\v{e}-Galton's basic fingerprint
patterns as models to be mathematically represented by the phase portrait of
an autonomous planar system of ordinary differential equations, which in an
equivalent form can be derived from $(1)$ and written as%
\begin{equation}
\overset{\cdot }{x}=f(x)  \label{2}
\end{equation}%
According to Newton's notation, the dot above $x$ means derivative with
respect to the independent variable (time, in general) and $f$ is a function
which the class of regularity is at least $C^{1}$ on $%
\mathbb{R}
^{2}.$ We know that, up to an appropriate rescaling of time, system $(2)$
defines a dynamical system on the plane, that is for which solutions are
uniquely defined for \textit{all} times. Practically, only a restriction of
it will be however considered, namely a restriction to a compact,
simply-connected subset $U$ of the plane which will be identified to the
cylindrical projection upon that plane of the inner surface of the last
phalanx of the finger in question (cf. chapter III of Galton's book on
methods of printing). For preliminary modelling, by mathematical
representation we roughly mean existence of a sufficiently smooth, not
necessarily computable, one-to-one transformation,\footnote{The term \textit{smooth deformation} of phase portraits will be used in the
remaining part of the paper. Continuity of the transformation and its
inverse is the minimum required. However, as we shall see, to preserve the
nature (but not necessarily the position) of each singular point, the degree
of regularity has to be increased when a classification up to homeomorphism
(topological equivalence) cannot distinguish between two singular points, as
for example a Node and a Focus.}\ by means of which trajectories of $(2)$ on
the \textquotedblleft restricted phase space\textquotedblright\ $U$ could be
mapped\ onto a simplified version of the fingerprint's streams of ridges,
i.e. the emerging fingerprint's general orientation feature and
corresponding \textquotedblleft trajectories\textquotedblright\ obtained,
for example, via a gradient-based or, in a more sophisticated way, a
filter-based method. In fact, a large literature is available on the subject
of computer-based methods for fingerprints features extraction,\footnote{%
We don't give a survey as it would load these pages too heavily to present
such technical methods here.} but for our sense of vision, whose accuracy in
pattern recognition cannot in principle be matched by any computer-based
method, Galton already pointed out that when contemplating a fingerprint,
the (unaided) eye is guided merely by the general appearance, while actually
the object under study can be much more complex. Moreover, what may still
bias the study is that \textit{a complex pattern is capable of suggesting
various readings, as the figuring on a wall-paper my suggest a variety of
forms and faces to those who have such fancies}. A simplified version,
whatever it may be, remains in fact a purely subjective notion. This,
however, does not prevent it from being of a certain utility for primary
classification purposes, and as we shall see, even in simplified form, and
no regard being paid to the fact whether the wanted transformation is
preserving or not orientation along trajectories, it will make it certain in
most cases that finality will never be perfectly reached, as the following
subsection shows.

\subsection{Searching for the Delta}

\begin{figure}
\begin{center}
\includegraphics{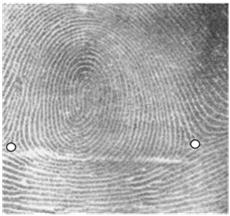} 
\end{center}
\caption{\it A Delta at both sides of
the Whorl.}
\label{delta}
\end{figure}

To model the last six categories of fingerprints according to Purkyn\v{e}'s
standards, and besides core points which will be discussed later, a quite
obvious (or at least what it seems to be) singular point is the Delta, as
called in Henry's classification and appears in Fig.\ref{delta}, such a
point being a basic singular point in impressions of the Loop and Whorl
types. It may be formed \textit{either by the bifurcation of a single ridge,
or by the abrupt divergence of two ridges that hitherto had run side by
side. }And as pointed out by Galton, following Purkyn\v{e}, in the Spiral,
the Ellipse, the Circle and the Double Whorl, \textit{\textquotedblleft
triangles\textquotedblright\ may be seen at the points where the divergence
begins between the transverse and the arched lines, and at both sides. }In
the language of dynamical systems, a Delta can be roughly seen as a singular
point with three hyperbolic sectors. So, the question is the following: is
there any singular point with three hyperbolic sectors?

\bigskip

Of course, in the first species class of singularities, there is no such
point. Note that a Saddle is a point with (a deleted neighborhood consisting
of) four hyperbolic sectors. As for singular points of the second species,
at least those from Fig.\ref{singular}, none looks like a Delta. In fact, as
said before, without searching to draw up an exhaustive list, such a point
cannot exist as it follows immediately from Bendixson formula above:

\begin{proposition}
\textit{If a singular point of system }$(2)$\textit{\ has three (and in
general, an odd number of) sectors, they can not all be of the hyperbolic
type.}
\end{proposition}

\bigskip

To take a closer look, the curious reader can suppose given a singular point
with three hyperbolic sectors and a defined direction of the flow for one of
the three sectors. He will then necessarily find himself faced with
incompatibly oriented trajectories, that is a situation where opposing
tangent vectors appear in an infinitely small region of the phase space.%
\footnote{%
At the beginning of chapter IV on ridges and their use, Galton talks about
minute ridges that bear a \textit{superficial} resemblance to those made on
sand by wind or flowing water. Indeed, it is a deceptive resemblance as it
is known for example in fluid dynamics, and unless there is an obstacle,
there can be no delta-like motion in a natural occurring flow.}
Mathematically, denoting by $S_{i},$ $i=1,2,3,$ the three separatrices, there
exist $%
i\in \{1,2,3\},\ y\in S_{i}$ and $\varepsilon >0$ such that, for any neighborhood 
$U$ of $y,$ there exist $\delta >0$ and $z\in U$ satisfying%
\begin{equation*}
\left\Vert \overset{\cdot }{x}(t,y)-\overset{\cdot }{x}(t,z)\right\Vert >\varepsilon \text{ \ for all \ }t\in
]-\delta ,\delta \lbrack
\end{equation*}%
$x(t,y)$ being the solution of the initial value problem%
\begin{equation*}
\left\{ 
\begin{array}{l}
\overset{\cdot }{x}=f(x) \\ 
x(0)=y%
\end{array}%
\right.
\end{equation*}%
and $\left\Vert \cdot \right\Vert $ the Euclidean norm on $%
\mathbb{R}
^{2}.$ This contradicts the regularity of the vector field $f$ and then the fundamental theorem on dependence on initial
conditions (and parameters, eventually) which, being based on Gronwall's
lemma, states that this dependence is - roughly speaking - as continuous as
the function $f.$ The general result for an arbitrary odd number $n>3$ of
hyperbolic sectors holds by the same reasoning.

\bigskip

The fact that there is no three-hyperbolic-sector point in smooth dynamical
systems is at the same time frustrating and fascinating. Frustrating because
the modelling process, still at an embryonic stage of realisation, seems to
be blocked. Fascinating when you realise that life has a lot more
imagination in printing than a differential equation can do!

\subsection[The Cusp, a faulty point but the best available]{The Cusp, a faulty point but the best available\protect\footnote{%
The title is inspired by the paragraph \textit{\textquotedblleft Fraternity,
a faulty word but the best available\textquotedblright }, from Galton's
chapter XI on heredity.}}

\begin{figure}
\begin{center}
\includegraphics[scale=0.7]{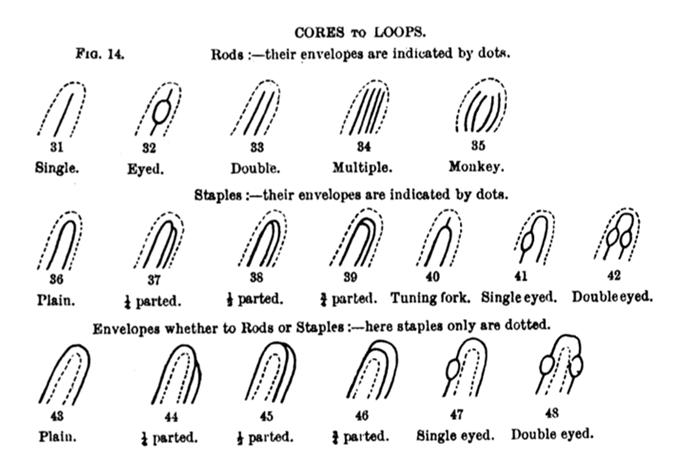} 
\end{center}
\caption{\it Cores to Loops, as captured from Plate 8 of Galton's book. At top left, the Single.}
\label{cores}
\end{figure}

As it has been mathematically shown, no singular point can match the Delta.
To overcome the problem, a well-known clever trick consists in doubling the
fingerprint's orientation field, thus transforming a Delta into a Saddle and
vice versa, after reconstruction. The same idea works for a special type of
cores, namely what is called \textquotedblleft Single\textquotedblright\ by
Galton, transforming it into a Focus, which is undoubtedly one of the best
points to be securely distributed by two Saddles, i.e. without changing its
nature. However, a Focus is not always to be expected from an orientation
doubling if the rich variety of cores is taken into account (see Fig.\ref%
{cores}), thus showing the limitation of this technique. Trying to see if
the task can otherwise be achieved, and as a three-sector singular point, on
could have thought for instance of the Saddle-Node as an approximation;
however, as can be easily seen from the structure of the Spiral or the
Elliptical/Circular Whorl, the parabolic sector of the Saddle-Node is not
appropriate to approach the center of the Whorl.\footnote{%
According to Definition 4.1 (cf. paragraph 4.4), I am currently trying,
however, to identify a family of Saddle-Node points arising in the
codimension-two Bogdanov-Takens bifurcation (Fig.\ref{BT}, $\mu _{1}=0$)
with some parted cores (Fig.\ref{cores}, 44-46).} In fact, in this quest for
the Delta, it becomes more and more certain that the goal will never be
reached by the path hitherto pursued, namely, to seek at all costs to
identify the Delta (and any pattern in general) to the \textquotedblleft
whole\textquotedblright\ of a singular point. I then had the idea of
{}{}considering only a \textquotedblleft portion\textquotedblright\ of it,
and I almost immediately realised that, for example, from a well-chosen cut
of the Cusp, it emerges something close to a Delta,\footnote{%
Such a cut can be justified only if one agrees that the hidden area can be
easily guessed and topologically reconstructed for all types of fingerprints.%
} as appears in Fig.\ref{circular}-\ref{spiral}. Then, as a preliminary
modelling of the corresponding fingerprint, one can place two Cusps
\textquotedblleft face to face\textquotedblright\ on the boundary of a
finger-shaped domain and put in the middle a Center for the
Elliptical/Circular Whorl, and a Focus for the Spiral. The result is the
following:

\subsection{Preliminary modelling of the Circular/Elliptical Whorl, the
Spiral, the Arch and the Loop}

\subsubsection*{The Circular/Elliptical Whorl and the Spiral}

\begin{figure}
\begin{center}
\includegraphics[scale=0.7]{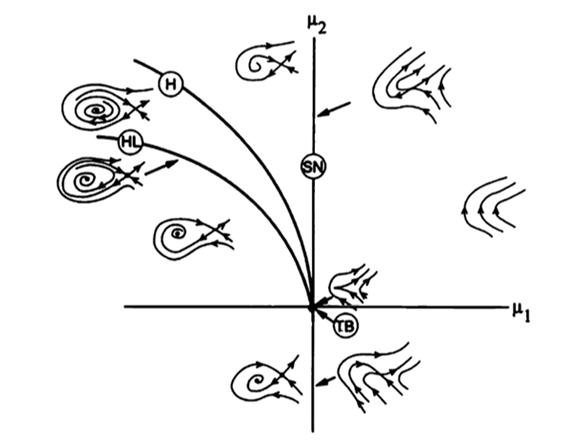} 
\end{center}
\caption{\it As captured from Perko's
book: bifurcation set and corresponding phase portraits for the
system $\protect\overset{\cdot \cdot }{x}=\protect\mu _{1}+\protect\mu _{2}%
\protect\overset{\cdot }{x}+x^{2}+x\protect\overset{\cdot }{x},$ \textit{TB,
H, SN and HL being for Takens-Bogdanov, Hopf, Saddle-Node and
Homoclinic-Loop bifurcations, respectively}.}
\label{BT}
\end{figure}


\begin{figure}
\begin{center}
\includegraphics[scale=1.2]{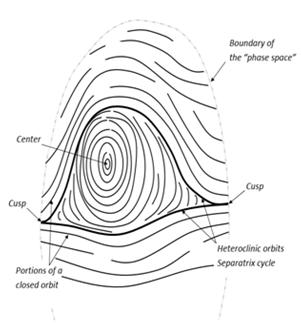} 
\end{center}
\caption{\it Preliminary modelling of the Circular/Elliptical Whorl.}
\label{circular}
\end{figure}


\begin{figure}
\begin{center}
\includegraphics{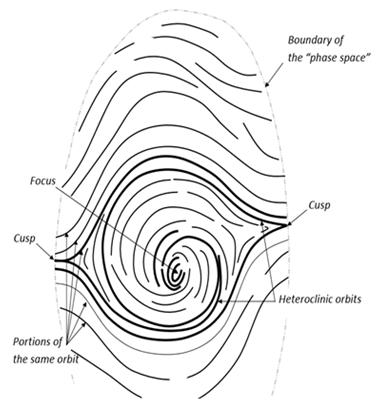} 
\end{center}
\caption{\it Preliminary modelling of the Spiral.}
\label{spiral}
\end{figure}


In Fig.\ref{circular} (in which, as well as in all sketched phase portraits,
and to better imitate a fingerprint, solution curves are deliberately
represented by irregular discontinuous lines), we have two (partial) local
phase portraits corresponding to a Cusp, and a local phase portrait
corresponding to a Center. These three local phase portraits have to be
connected to obtain a global phase portrait, whence the notion of connexion
mentioned in the introduction. Intuitively, this fact can be achieved only
if the curves on either side of the Center meet symmetrically one-to-one to
form closed orbits, the separatrices acting as heteroclinic orbits joining
the Cusps and forming what is called a separatrix cycle. We obtain in fact a
similar configuration as for the stable equilibrium of the undamped
pendulum, where the Center is now served by Cusps instead of Saddles. As for
the Spiral in Fig.\ref{spiral}, a possible configuration is simply that of
Fig.\ref{basic}, but with Cusps instead of Saddles, the origin being a
global attractor\textit{\ }for all trajectories, except those corresponding
to the stable manifolds of the Cusps. Another modelling phase portrait
which, up to a smooth deformation, and if we exclude phase portraits with
(bifurcating) limit cycles, seems to be the only other possibility for a
configuration where a Focus is trapped between two symmetric Cusps, is that
for which the separatrix cycle is approached by the unbounded curves and the
inner spiral as a limit set (the set of cluster points of the
forward/backward orbit), no regard being paid to the direction of the
\textquotedblleft twist\textquotedblright\ or to the number of
\textquotedblleft turns\textquotedblright . Theoretically, the state of the
limit set is reached after an infinite number of forwards or backwards
turns. As the case may be, we will have a stable, unstable or semistable
separatrix cycle. It may be asked why the well-known phase portraits
corresponding to an undamped pendulum or to a simple case (zero drive
strength) of the damped pendulum have not been directly considered as
preliminary models for the Elliptical/Circular Whorl or the Spiral,
respectively. The reason in fact why I did not adopt such a configuration is
that, unlike a Cusp, Saddle' separatrices do not meet tangentially as
required by the enveloping ridges of a pattern, a behaviour that can be
easily observed from a close-up of the Delta in Fig.\ref{delta}. Besides,
considering a Saddle instead of a Cusp supposes neglecting a whole
hyperbolic sector, whereas the two sectors clearly appear in the case of a
Cusp, only small portions of (eventually closed) curves have been ignored.
Letting Galton comment on these figures, he would probably say: \textit{What
seemed before to be a vague and bewildering maze of lineations over which
the glance wandered distractedly, seeking in vain for a point on which to
fix itself, now suddenly assumes the shape of a sharply-defined phase
portrait.\footnote{%
In fact, in [1, p.69], Galton was speaking about the \textit{%
outlines} of patterns and how they can be accurately drawn. I just took the
liberty to replace \textquotedblleft figure\textquotedblright\ with
\textquotedblleft phase portrait\textquotedblright\ !}}

\bigskip

The question that remains however is the following: Given in general
fragments of phase portraits, how to rigorously carry out (all possible
cases of) a connexion? And once it has been accomplished, how to validate\
it, that is to construct algebraically a dynamical system whose phase
portrait is a smooth deformation of the connexion? This is the question to
which we try to provide some elements of answer in section 5. In the
following, we resume the preliminary modelling for the simple cases in
Arches and Loops, namely the Transverse flexures, the Central Longitudinal
Stria and the Oblique Sinus, where the Cusp will be found to give
considerable help. The remaining cases of the Oblique Stria, the Almond, the
Composite Spiral (in a sense to be specified) and the Double Whorl will be
discussed in the next subsection.

\subsubsection*{The Arch}

\begin{figure}
\begin{center}
\includegraphics{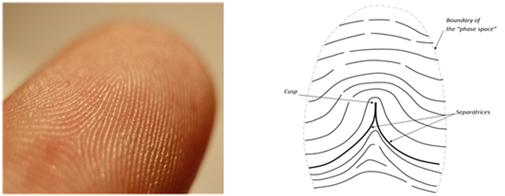} 
\end{center}
\caption{\it Tented Arch on male
finger (Wikimedia Commons) and associated phase portrait.}
\label{arch}
\end{figure}

For fingerprints of the Arch type, especially those with the transverse
flexures for example, \textit{where ridges are arranged transversely in
beautiful order}, and as a preliminary modelling, a straightforward
topological similitude can be made with straight curves running transversely
from one side of the phalanx to the other. In terms of the well-known
rectification theorem from dynamical systems, the topologically equivalent
phase portrait can be seen as a magnification of the flow near a nonsingular
point $x$, after rectification, that is after having applied \textit{a
change of coordinates for a region around }$x$\textit{\ where the vector
field }$f$\textit{\ becomes a series of parallel vectors of the same
magnitude}. In other words, a Plain Arch can be merely seen as a smooth
deformation of a line. However, at another level, a pattern of Arch type can
be a bit more complex, as for example it is the case for the Central
Longitudinal Stria, where the configuration is nearly the same as in the
previous case, but with a small difference: \textit{a perpendicular stria is
enclosed within the transverse furrows, as if it were a nucleus}. In
Galton's language, this case, extremely rare on the thumb,\footnote{\textit{I
do not remember ever to have seen it there}, says Galton in [1, p.75].} can be included within the \textit{Tented Arch} for which an
approximate, but apparently fairly correct, phase portrait can be expressed
by a Cusp (Fig.\ref{arch}).

\subsubsection*{The Loop}

\begin{figure}
\begin{center}
\includegraphics{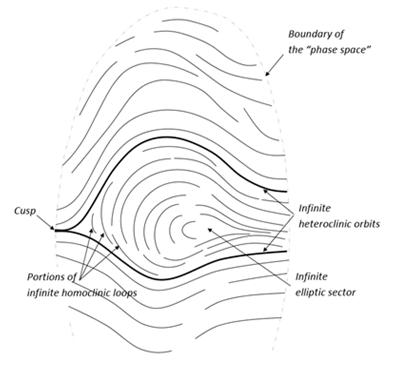} 
\end{center}
\caption{\it Preliminary modelling of
the Loop.}
\label{Loop}
\end{figure}

For the Oblique Sinus, where \textit{an oblique line recurves towards the
side from which it started, accompanied by several others, all recurved in
the same way,} and as a preliminary modelling within the principal head of
Loops, a simplified version can be seen as a smooth deformation of a phase
portrait where a Cusp cohabits with an infinite singular point\footnote{%
With all due respect to Hermite who, unlike Poincar\'{e}, opposes a
vocabulary which he finds too colourful, like \textit{points at infinity} in
projective geometry!} with (at least) an elliptic sector. To be placed
along the vertex of the oblique sinus of the furrows distribution, the
specific need of an infinite point comes from the fact that the unbounded
Cusp' separatrices have to be flexed in order to cover the recurve lines as
enveloping ridges. And for the elliptic sector, as the name suggests, it
seems that no solution curve can simulate the shape of a loop better than
that of an elliptic sector. The general feature of a Loop's phase portrait
is sketched in Fig.\ref{Loop}.

The main complaint one might have concerning this model, however, is the
lack of the core point. For a more accurate phase portrait, in the case of a
Single core, another adjustable slightly sloping Cusp can be integrated
within the picture. Otherwise, getting even more accurate, but requiring a
larger phase space, an interesting idea to be elaborated in future work is
proposed in subsection 4.4.3, where the same problem is encountered in the
Double Whorl modelling.

\subsection{Preliminary modelling of the Oblique Stria, the Almond and the
Composite Spiral}

\subsubsection*{The Oblique Stria}

The pattern is defined as Transverse flexures where \textit{a solitary line
runs from one or other of the two sides of the finger, passing obliquely
between the transverse curves, and ending near the middle}. In the language
of dynamical systems, this can be achieved by considering (a smooth
deformation of) the phase portrait associated to the Tented Arch, but with a
deleted separatrix for the Cusp. Doing so, and for reasons that will become
clear later, we introduce a modified definition of what will be meant from
now on by a phase portrait, whose a smooth deformation is to be associated
to the ridge flow of a fingerprint:

\bigskip

\begin{definition}
Given a compact, simply-connected domain $U\subset 
\mathbb{R}
^{2},$ and a nonempty finite subset $V$ of $U,$ a phase portrait on $U$ with
respect to $V$ is defined to be the set%
\begin{equation*}
P=\left\{ \varphi _{\pm t}(x),\ t\geq 0,\ x\in V\right\} \cap U
\end{equation*}%
where $\varphi _{t}$ is the (one-parameter) flow generated by a vector field 
$f\in C^{1}(%
\mathbb{R}
^{2}),$ and the symbol $\pm $ means that, for an initial condition $x$ and
positive times $t,$ one has the choice to consider the solution curve either
of the system $\overset{\cdot }{x}=f(x)$ or $\overset{\cdot }{x}=-f(x).$
\end{definition}

\bigskip

Remember, $U$ is practically identified to the cylindrical projection upon a
plane of the inner surface of the last phalanx of the finger at issue; it
has been roughly considered so far as a restricted phase space in the
preliminary modelling of fingerprints. As for the underlying vector field,
it should be stressed that $f$ is defined on the whole plane, the modelling
phase portrait being created by the flow generated by the restriction of $f$
to $U,$ according to Definition 4.1. In this regards, it does not matter how $%
f $ behaves outside $U$ as long as it works inside, i.e. it gives the
desired phase portrait on $U.$

\bigskip

Now, according to Definition 4.1, any initial condition can be chosen to
generate an Oblique Stria phase portrait, except a point located on the
Cusp' separatrix to be ignored. Especially, denoting by $W^{s}$ (resp. $%
W^{u})$ the stable (resp. unstable) manifold of the Cusp and by $\partial U$ the boundary of $U$, if we let $x_{0}$ be a
nonsingular point so that%
\begin{equation*}
W^{s}\cap \partial U=\left\{ x_{0}\right\} \ \ \text{(resp. }W^{u}\cap
\partial U=\left\{ x_{0}\right\} )\ ,
\end{equation*}%
the trajectory%
\begin{equation*}
T^{+}=\left\{ \varphi _{t}(x_{0}),\ t\geq 0\right\} \ \ \text{(resp. }%
T^{-}=\left\{ \varphi _{t}(x_{0}),\ t\leq 0\right\}
\end{equation*}%
will stand for a smooth deformation of the oblique stria (the ridge, not the pattern), depending on which
side of the finger the solitary line is running from.

\subsubsection*{The Almond}

What is called Almond by Purkyn\v{e} and described as \textit{an Oblique
Sinus enclosing an almond-shaped figure, blunt above, pointed below, and
formed of concentric furrows}, is in fact a compound pattern called Circlet
in Loop by Galton, or, as the case may be, Spiral in Loop. As pointed out by
Galton, Whorls enclosed within Loops are by far the commonest pattern among
the compound category. So, trying to encrust a whorled pattern into a Loop,
remember the interesting cases we have met in the codimension-two
Bogdanov-Takens bifurcation (Fig.\ref{BT}, $\mu _{1}<0)$, where a Focus is
connected to a Saddle, directly by a heteroclinic orbit (Fig.\ref%
{Spiral-Loop}), via a homoclinic loop or a limit cycle, a configuration
which can be used for the Spiral in Loop. The case of the Circlet in Loop
can be approached in the same manner where a Center is connected to a Saddle
via a homoclinic loop (Fig.\ref{Circlet-Loop}). To complete the picture, a
Cusp has to be added, where the lower separatrix has to be connected to an
infinite Node in order to delimit the pattern, and according to Definition
4.1, a Saddle' separatrix has to be deleted, when some trajectories (within
the same phase portrait) have to be drawn backwards.

\begin{figure}
\begin{center}
\includegraphics{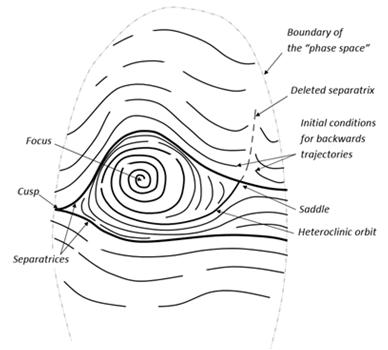} 
\end{center}
\caption{\it Preliminary modelling of the Spiral in Loop.}
\label{Spiral-Loop}
\end{figure}


\begin{figure}
\begin{center}
\includegraphics{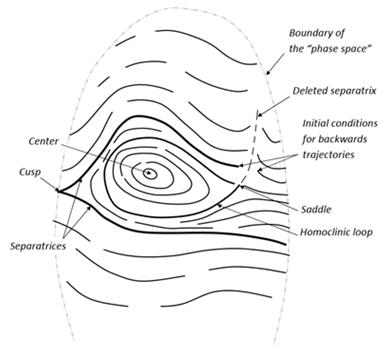} 
\end{center}
\caption{\it Preliminary modelling of the Circlet in Loop.}
\label{Circlet-Loop}
\end{figure}


\subsubsection*{The Composite Spiral}

\begin{figure}
\begin{center}
\includegraphics{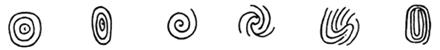} 
\end{center}
\caption{\it Cores to Whorls, from
left to right, in Galton's language: Circles, Ellipses, Spiral, Twist,
Plait, Deep Spiral \protect\cite{Galton}.}
\label{cores-whorls}
\end{figure}

Remember, the term \textquotedblleft Spiral\textquotedblright\ has been used
so far in the usual geometric sense. However, as pointed out by Galton, if
the term \textit{conveys a well-defined general idea, there are four
concrete forms of it which admit of being verbally distinguished}: the
(simple) Spiral, the Twist, the Plait and the Deep Spiral. In addition to
Circles and Ellipses, they appear as Cores to Whorls in Plate 8 - Fig.15 of
Galton's book (see Fig.\ref{cores-whorls}).

\bigskip

In Purkyn\v{e}'s Commentatio, only the Spiral and the Twist are considered
as spirals, the former being classified as simple and the latter as
composite. The Plait is called Double Whorl therein and no reference seems
to be made to the Deep Spiral. Trying to identify each of the
\textquotedblleft complex\textquotedblright\ spirals with a singular point,
I realised that they were ranked in ascending order of difficulty, as if
Galton had already tried to identify them with a phase portrait! In fact, I
could not do anything for the Deep Spiral, however, for the Twist, and to a
lesser extent the Plait, two local phase portraits have been found to be of
some help, as will be explained in the following.

\paragraph{The Twist}

\begin{figure}
\begin{center}
\includegraphics[scale=0.7]{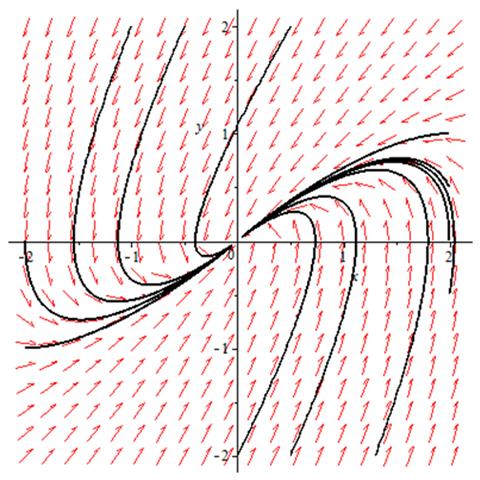} 
\end{center}
\caption{\it Computer-drawn phase
portrait of an Improper Node.}
\label{Twist-Maple}
\end{figure}

\begin{figure}
\begin{center}
\includegraphics{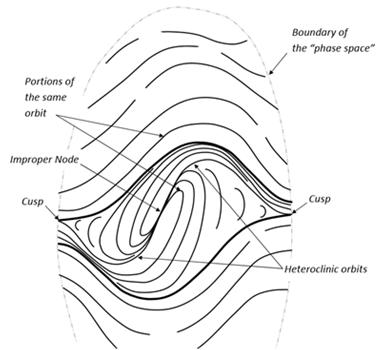} 
\end{center}
\caption{\it Preliminary modelling
of the Twist.}
\label{Twist}
\end{figure}

The spiriform of the pattern is described by Purkyn\v{e} as being \textit{%
made up of several lines proceeding from the same centre, or of lines
branching at intervals and twisted upon themselves}. The best I could do for
such a pattern is to draw the phase portrait of an Improper Node, a singular
point known for being a transitional case between a Node and a Focus. In
fact, algebraically, for linear differential planar systems (Fig.\ref%
{diagram}), the Improper Node is located on the parabola $\tau
^{2}-4\triangle =0,$ for which a double root is exhibited by the
characteristic polynomial, passing continuously from the case of two
distinct roots with the same sign to complex conjugate ones. And as said
before, this is the reason why a homeomorphism cannot distinguish between a
Node and a Focus. Geometrically, through the Maple-drawn phase portrait (Fig.%
\ref{Twist-Maple}) and the hand-drawn preliminary model (Fig.\ref{Twist}),
the judgement is left to the reader on whether an Improper Node can be a
good candidate for a Twist or not. Finally, to obtain the whole picture, we
don't forget to place (and connect) the Improper Node between two Cusps, as
said by Purkyn\v{e} when describing his Composite Spiral: \textit{At either
side, where the spiral is contiguous to the place at which the straight and
curved lines begin to diverge, in order to enclose it, two triangles are
formed, just like the single one that is formed at the side of the Oblique
Sinus}.

\paragraph{The Double Whorl, or the need to increase dimension}

Also a Plait (and sometimes an Overlap) for Galton, it is described as 
\textit{a curious effect where two systems of ridges that roll together, end
bluntly, the end of the one system running right into a hollow curve of the
other, and there stopping short; it seems, at the first glance, to run
beneath it, as if it were a plait. This mode of ending forms a singular
contrast to the Spiral and the Twist, where the ridges twist themselves into
a point}. And for Purkyn\v{e}, \textit{one portion of the transverse lines
runs forward with a bend and recurves upon itself with a half turn, and is
embraced by another portion which proceeds from the other side in the same
way. This produces a doubly twisted figure which is rarely met with except
on the thumb, fore, and ring fingers}.

\bigskip

\begin{figure}
\begin{center}
\includegraphics{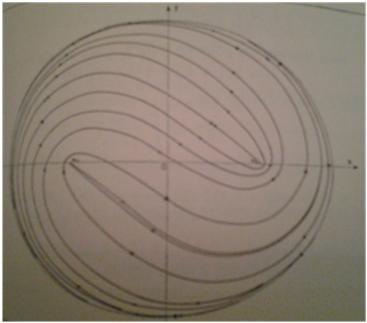} 
\end{center}
\caption{\it A plait-like behaviour
of trajectories on a Poincar\'{e} section from a special case of the
restricted planar three-body-problem.}
\label{Plait}
\end{figure}

\begin{figure}
\begin{center}
\includegraphics{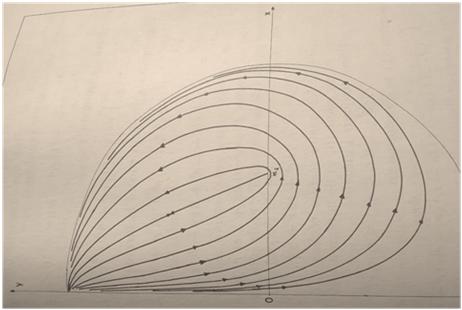} 
\end{center}
\caption{\it A loop-like behaviour of trajectories on a Poincar\'{e} section from a special case of the restricted planar three-body-problem. Compare this pattern to the one on C\'{e}sar Baldaccini's Thumb from the Theme of the Hand exhibition!}
\label{Loop1}
\end{figure}

After numberless observations, I failed completely in trying to translate
this in terms of solution curves of a planar differential system; but I have
been vaguely conscious of remembering that a similar figure can be obtained
on a Poincar\'{e} section when considering a two-degree-of-freedom
Hamiltonian, i.e. a four dimensional conservative system. Indeed, such a
type of behaviour can be encountered in classical mechanics when dealing,
for example, with the Lagrange problem of equal masses, a special case of
the well-known restricted planar three-body problem when Coriolis' force is
neglected. Without going into details, the Double Whorl, roughly expressed
in Fig.\ref{Plait}, is taken from an international postgraduate course given
by A. Deprit in August 1960 at the Universit\'{e} Libre de Bruxelles.%
\footnote{%
As a young researcher, I have met the late Prof. Andr\'{e} Deprit for the
first time in the summer of 1998 in Prague. The private communications that
I had the honor to tie with him, mainly on normal form theory and celestial
mechanics, have left their imprint until today.} The plait-like behaviour
is undergone by a special family of solution curves on the invariant $(x,y)$
plane for the system with the Lagrangian function%
\begin{equation*}
L=\frac{1}{2}(\overset{\cdot }{x}+\overset{\cdot }{y})+\frac{1}{4}%
(r_{13}^{2}+r_{23}^{2})+\frac{1}{2}(1/r_{13}+1/r_{23})
\end{equation*}%
$r_{13}$ being the distance of the moving particle to the origin and $r_{23}$
the distance between the two masses within the Keplerian motion. The \textit{%
mode of ending} perfectly described by Galton corresponds in fact to the
position of the equal masses, and the \textit{two systems of ridges that
roll together }are portrayed by a family of six asymptotic orbits, beginning
with an ejection orbit from one of the two masses, passing through orbits
that have missed the double-shock, thus continuing to revolve around that
same mass, ending finally with an ejection orbit from the other mass. At
another energy level, solution curves behaviour could be also of some help
for the preliminary modelling of the Loop, where the missing core in Fig.\ref%
{Loop} is now portrayed by a collision orbit with one of the two masses,
then reaching back to join asymptotically one of the five well-known
Lagrange points (Fig.\ref{Loop1}). However, to complete the picture, and if
all is confined in an ellipse with the equal masses positions as focal
points, I don't know if there exist off-ellipse initial conditions leading
to a (partial) delta-like motion, let alone how to proceed to a connexion.
Carrying out a connexion from given data in the general case is the main
subject of the following section, and as can be suggested for future work,
besides index theory on a two-dimensional surface and multivariate
interpolation theory, bifurcation theory is another framework within which
the problem could be formulated and dealt with. More precisely, the rather
fascinating question of studying all possible connections of local phase
portraits could be seen as a bifurcation subproblem, i.e. a global
bifurcation problem with the constraint that all involved singular points
are preset to predetermined geographical positions and behavioural natures.%
\footnote{%
Obviously, the term bifurcation, at least in its dynamical sense, appears
nowhere in Galton's book, but when the author talks about sets of concentric
circles or ellipses, pointing out in p. 77 that\textit{\ they are rarely so
in a strict sense throughout the pattern, usually breaking away into a more
or less spiriform arrangement}, we dare to wonder if he was already
anticipating the well-known Andronov-Hopf bifurcation! Better yet, the 
\textit{transitional cases} concisely described by Galton, as those between
a Tented Arch and a Loop, could be expressed in a dynamical context as a
bifurcation problem.}



\part{\protect\Large Connecting singular points within structure-preserving
normal forms}

\section{Some directions in connecting phase portraits}

Recapitulating, the way hitherto pursued in modelling consisted in locating
basic fingerprint patterns, identifying them to (portions of) local phase
portraits, then, guided by intuition,\footnote{%
By \textit{intuition}, we allude to all these rules from dynamical systems
to be respected when drawing up phase portraits, beginning with the
fundamental theorem of existence and unicity, and arriving at the most
elaborate theories, like Poincar\'{e}'s index, consequents, contacts,
centers and limit cycles theories [9-11]. For example, it would be inconceivable to draw a cycle without
any singularity inside, to delimit a Center by a limit cycle, or to try to
connect a Node and a Focus without resorting to other singular points.}
making a connexion to obtain a global phase portrait. Although this results
in interesting phase portraits, it must be recognised that the whole process
is a craftsman's task, and no systematic approach has been proposed.
Besides, there is no guarantee for the correspondence via a smooth
deformation between the object thus constructed and the phase portrait of an
explicitly written (in terms of elementary functions) dynamical system. As
for the problem of existence of such a phase portrait in the general case,
i.e. for any number (and nature!) of singular points, nothing is known so
far, especially in case we fail in the connexion process. In the next
paragraph, we deal with the problem within Poincar\'{e}'s index theory to
show that a connexion is generally impossible, but when it exists,
obviously, a connexion is not necessarily unique, as it will be shown
through a simple illustrating example. In this regards, some results from
index theory will be briefly recalled to deal with the question of existence
in the general case, and as it will be seen, the explicit construction of
the underlying dynamical system, i.e. the ultimate goal of associating a
(smooth deformation\footnote{%
A link to smooth deformations could be made from Galton's work, as for
example when he says in p. 75: \textit{Perhaps the best general rule in
selecting standard outlines, is to limit them to such as cannot be turned
into any other by viewing them in an altered aspect, as upside down or from
the back, or by magnifying or deforming them, whether it be through
stretching, shrinking, or puckering any part of them.}} of a) differential
system's phase portrait to a fingerprint, according to Definition 4.1, leads
to a dreadful and generally insoluble problem.

\subsection[Index theory revisited]{Index theory revisited\protect\footnote{%
As the reader will notice, the purpose for which the theory is recalled,
i.e. the problem of connecting local phase portraits, has nothing whatever
to do with the way Poincar\'{e}'s index is used for singular points
extraction from fingerprints coarse orientation fields.}}

\subsubsection{Towards Poincar\'{e} index theorem}

To give a complete account of the theory would require a chapter, so I have
opted for a sketch. Following Poincar\'{e}, we recall the definition of the
index $I_{f}(x_{0})$ of an isolated singular point of a $C^{1}$ vector field 
$f,$ defined on an open subset $U$ of $%
\mathbb{R}
^{2},$ as being the index $I_{f}(C)$ of any Jordan curve $C\subset U,$
containing $x_{0}$ and no other singular point of $f$ on its interior, which
is given by%
\begin{equation*}
I_{f}(x_{0})=I_{f}(C)=\frac{\triangle \theta }{2\pi }
\end{equation*}%
$\triangle \theta $ being the total change in the angle $\theta $ that the
vector $f=(f_{1},f_{2})^{T}$ makes with respect to the $x$-axis, i.e. the
change in%
\begin{equation*}
\theta (x,y)=\arctan \frac{f_{1}(x,y)}{f_{2}(x,y)}
\end{equation*}%
as the point $(x,y)$ traverses $C$ exactly once in the positive direction.
Explicitly, this can be computed as%
\begin{equation*}
I_{f}(x_{0})=\frac{1}{2\pi }\oint_{C}\frac{f_{1}df_{2}-f_{2}df_{1}}{%
f_{1}^{2}+f_{2}^{2}}
\end{equation*}%
In index language, with the same notation, the previously mentioned
Bendixson result (on the number of elliptic and hyperbolic sectors of a
singular point) can be stated for analytic $f$ as follows:

\begin{theorem}
\textbf{(Bendixson)}%
\begin{equation*}
I_{f}(x_{0})=1+\frac{e-h}{2}
\end{equation*}
\end{theorem}

\bigskip

It follows that, for instance, in the nondegenerate case, $I_{f}(x_{0})$ is $%
-1$ or $+1$ according to whether the singular point is or is not a
(topological) Saddle. Note also that the index of a Saddle-Node is zero.

\bigskip

With respect of the vector field $f,$ the index theory is extended to a
two-dimensional surface $S$ (i.e. a compact, two-dimensional, differentiable
manifold of class $C^{2},$ or nappe as would say Poincar\'{e} in French$),$
on which $f$ is supposed to have a finite number of singular points $%
x_{1},...,x_{m}.$ $I_{f}(S)$ is then defined as the sum of the indices at
each of the singular points:%
\begin{equation*}
I_{f}(S)=\sum_{i=1}^{n}I_{f}(x_{i})
\end{equation*}%
$I_{f}(x_{i})$ being defined relatively to the restriction of $f$ to some
chart, into the detail of which it is unnecessary here to enter.

\bigskip

As shown by Poincar\'{e} \cite{Poincare2} in chapter XVIII, following a work
in two parts of Kronecker \cite{Kronecker}, it is one of the most
interesting facts of the index theory that $I_{f}(S)$ is independent of the
vector field $f$ and only depends on the topology of the surface:

\begin{theorem}
\textbf{(Poincar\'{e} Index Theorem)}%
\begin{equation*}
I_{f}(S)=\chi (S)
\end{equation*}%
\textit{where }$\chi (S)=T+v-l$\textit{\ is the Euler characteristic
associated to a decomposition of }$S$\textit{\ into a number }$T$\textit{\
of curvilinear triangles, with a number }$v$\textit{\ of vertices and a
number }$l$\textit{\ of edges.}
\end{theorem}

\textit{\bigskip }

It can also be shown that%
\begin{equation*}
\chi (S)=2(1-p)
\end{equation*}%
where $p$ is the genus of $S$ (i.e. the maximum number of nonintersecting
closed curves than can be drawn on $S$ without dividing it into two separate
regions), thus leading to the topological invariance of $\chi (S)$ and then
of $I_{f}(S).$ As examples we have $I(S^{2})=2$ for the sphere, $I(T^{2})=0$
for the two-dimensional torus, $I(P)=1$ for the projective plane or $I(K)=0$
for the Klein bottle.

\bigskip

An immediate consequence of the Poincar\'{e} Index Theorem is the following:

\begin{corollary}
\textbf{(Poincar\'{e}) \ }\textit{Suppose that }$f$\textit{\ is an analytic
vector field on an analytic, two-dimensional surface }$S$\textit{\ of genus }%
$p$\textit{\ and that }$f$\textit{\ has only hyperbolic singular points,
i.e. isolated Saddles, Nodes and Foci, on }$S.$\textit{\ Then}%
\begin{equation*}
n+f-s=2(1-p)
\end{equation*}%
\textit{where }$n,$\textit{\ }$f$\textit{\ and }$s$\textit{\ are the number
of Nodes, Foci and Saddles on }$S$\textit{\ respectively.}
\end{corollary}

\subsubsection{Link to the connexion problem}

\begin{figure}
\begin{center}
\includegraphics[scale=0.7]{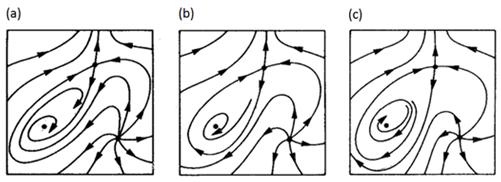} 
\end{center}
\caption{\it Qualitatively different connexions
for the same singularities.$^{19}$}
\label{connexion}
\end{figure}

Now, as an answer to the existence problem of a connexion, at least in the
hyperbolic case, it can be easily seen that for a given surface and a set of
Nodes, Foci and Saddles, which the respective numbers do not satisfy
condition of Corollary 5, it is impossible to carry out a connexion. For the
nonhyperbolic, eventually degenerate, case, any collection of singular
points whose the sum of indices does not match the Euler characteristic of
the underlying surface can not be connected. In other terms, a connexion of
index-theory-incompatible singular points with respect to a surface is
simply a faulty phase portrait on that surface. As for uniqueness, we simply
consider a portion of a phase portrait made up of local phase portraits from
a Saddle, a Node and a Focus, but in three qualitatively different ways (Fig.%
\ref{connexion}).

\bigskip

In (a), the Focus is joined by heteroclinic orbits from the Node and the
Saddle, the laters being themselves joined by two heteroclinic orbits
corresponding to the stable manifolds of the Saddle; in (b), the Focus is
enclosed within an homoclinic orbit corresponding to a separatrix cycle; and
in (c), the Focus is surrounded by an unstable limit cycle. The three phase
portraits, obviously, are not topologically equivalent.{\let\thefootnote\relax\footnote{$^{19}$These phase portraits are taken from an anonymous pdf document; I apologise
to the unknown author(s) for using these drawings to illustrate the
non-uniqueness of a connexion.}}

\bigskip

As for the general existence problem, it seems mathematically possible, but
highly unlikely, that a (large) number of pre-assigned singularities,%
\footnote{%
By pre-assigned singularities, we mean data singular points which are
intrinsically considered as such, i.e. as independent geometrical objects,
or as local phase portraits, like those in Fig.\ref{singular}, no regard
being paid for the moment to the explicit form of vector fields exhibiting
such singularities.} randomly distributed on a surface, yet
index-theory-compatible, can be connected. In other terms, we have a
necessary condition from Poincar\'{e} Index Theorem which seems to be
sufficient (to build a global phase portrait) only up to a redistribution of
the singular points on the surface. The second part of the problem is, once
a connexion is proved possible from the given data, find an
explicitly written differential system whose phase portrait is a smooth
deformation of one of the possible connections.\textit{\ }Such a reverse
problem\textit{\ }from index theory will be explicitly stated in the next
subsection as an interpolation problem.\ Consequently, whenever possible,
fingerprints discussed in section 4 will be assigned a simple planar
differential system as a primary model within the ALW classification.

\subsection{A bivariate interpolation problem}

\subsubsection{Robustness of singular points against addition of nonlinear
terms}

Let there be an open subset $U$ of the plane in which a finite number $n$ of
pre-assigned singular points $x_{i}$ have been successfully connected. The
problem of finding an explicit dynamical system%
\begin{equation*}
\overset{\cdot }{x}=f(x),\ x\in 
\mathbb{R}
^{2}
\end{equation*}%
with the same singularities on $U$ and whose phase portrait is a smooth
deformation of the connection in the vicinity of each singular point (but
not necessarily the case elsewhere), amounts to find a closed-form
sufficiently smooth (say, analytic) vector valued function $f$ whose zeros
in $U$ are exactly the $x_{i}$'s and Taylor polynomial at $x_{i},$ given
with a sufficiently high degree, defines a (local) differential system whose
trajectories exhibit the same behaviour as already pre-assigned for $x_{i},$
for all $i=1,...,n.$ By sufficiently high degree, we mean the smallest order
at which $f$ should be expanded so that the nature of the singular point
would be completely decided, and no additional high-order terms can destroy
this nature. In fact, as well-known within structural stability, under the
addition of smooth nonlinear terms to a linear system with a (global)
Center, the singularity may become a Focus or even a Center-Focus (in the
nonanalytic case). Besides, a homeomorphism does not distinguish between a
Node and a Focus, where a Saddle remains a (topological) Saddle under
continuously differentiable perturbations, as follows from the
Hartman-Grobman Theorem (1959). However, it follows from a result shown in
Coddington and Levinson \cite{CL} (resp. a theorem by Hartman (1960)) that a
Focus (resp. a Node) remains diffeomorphically unchanged under the addition
of (resp. twice) continuously differentiable nonlinear terms. This, indeed,
is the reason why the approach by Li and Yau \cite{LY} works, cores and
Deltas being respectively modelled by Foci and Saddles. As for the Center, a
large part has been devoted to the subject by Poincar\'{e} \cite{Poincare2},
chap. XI, where necessary and sufficient (yet difficult to implement)
conditions were derived to ensure the preservation of the center-like
nature. Following Poincar\'{e}, the simplest sufficient condition for a
nonlinear planar system to preserve the center nature is a symmetry with
respect to one or both of the axes, and this is probably the reason why the
approach by Ford \cite{Ford} works for typical flow configurations in fluid
dynamics. Therefore, recapitulating, if we restrict ourselves to the
nondegenerate case, given a distribution of Centers $C_{i},$ Foci $F_{i},$
Nodes $N_{i}$ and Saddles $S_{i}$ in $U$, it is clear that if there exists a
symmetric and twice continuously differentiable vector field $f$, vanishing
in $U$ only at these points, and for which $C_{i}$ (resp. $F_{i},$ $N_{i},$ $%
S_{i}$) is a Center (resp. Focus, Node, Saddle) of the linearised vector
field (i.e. defined by Taylor polynomial of degree 1), then the phase
portrait of the flow generated by $f$ stands for a connexion of the
pre-assigned data of singular points. A simple example in the nondegenerate
case is given below.

\subsubsection{Poisedness of the problem}

Let begin by considering the following actual situation: given a
configuration in the plane $(x,y)$ where a Center is placed at the origin
and two Saddles at the opposite points $\pm (\pi ,0),$ find a connexion in
the square $U=\{(x,y)\in 
\mathbb{R}
^{2}~/~\left\vert x\right\vert $ and $\left\vert y\right\vert \leq 4\}.$
Mathematically, this amounts to find a symmetric system,%
\begin{equation*}
\left\{ 
\begin{array}{ll}
\overset{\cdot }{x}= & f(x,y) \\ 
\overset{\cdot }{y}= & g(x,y)%
\end{array}%
\right.
\end{equation*}%
that is a system which is invariant under the transformation $%
(t,y)\longrightarrow (-t,-x)$ and/or $(t,y)\longrightarrow (-t,-y),$ where $%
f $ and $g$ are $C^{2}$-functions such that 
\begin{equation}
f(0,0)=f(\pm \pi ,0)=g(0,0)=g(\pm \pi ,0)=0\text{ \ }  \label{3}
\end{equation}%
and there is no other points of $U$ in which $f$ and $g$ vanish
simultaneously. Moreover, (cf. diagram of bifurcation, Fig.\ref{diagram}), 
\begin{equation}
\det (J_{\pm \pi })<0,~\text{{\small tr}}(J_{0})=0,\ \det (J_{0})>0,
\label{4}
\end{equation}%
$J_{\pm \pi }$ and $J_{0}$ being the Jacobian matrices at $\pm (\pi ,0)$ and
the origin, respectively.

\bigskip

As can be seen, the problem can be solved for $f(x,y)=y$ and $g(x,y)=-\sin
x. $ It can also be shown in fact that the phase portrait of the flow
generated by such a vector field is - up to a smooth deformation - the
unique connexion for a Center and two Saddles. However, as the reader has
certainly noticed, we purposely considered a problem of which we already
know the solution, namely the phase portrait of the undamped pendulum. But
in a more general context, condition (3) can be seen as a bivariate Lagrange
interpolation problem which is poised in the Haar space $C^{1}(U),$ i.e.
(allowing some flexibility in the definition) it can be solved in $C^{1}(U)$
for any given data of isolated points in $U,$ but not in the subspace of $%
C^{2}(U)$ of symmetric vector fields. As for conditions (3-4), they can be
seen as an advanced bivariate Hermite interpolation problem which is not
poised in $C^{1}(U),$ as follows from index theory. It should be stressed
that, for both problems, however, we are not necessarily dealing with
polynomial interpolation on the one hand, and on the other, no values have
been assigned to the partial derivatives of $f$ and $g$ at the singular
points, only a special distribution in the trace-determinant plane of the
Jacobian matrices is required.

\bigskip

The situation is dramatically different in the degenerate case, where
high-order derivatives are required to decide on the nature of the singular
point. And it is even more difficult to formulate such a problem in the
general case where a mixture of degenerate and nondegenerate points are
given as connected pre-assigned data. As will be shown, however, the problem
can be solved in some special situations where a restrictive number of
singularities are considered, as for the configurations \textquotedblleft
one Center/Focus/Improper Node - two Cusps\textquotedblright , but no simple
form has yet been found for the Loop (due to the interpolation point at
infinity) and the Double Whorl (due to the narrowness of the phase plane to
harbour such a behaviour).

\subsection{Introducing normal forms for fingerprints}

\subsubsection{A normal form for nilpotent planar systems}

A starting point to solve the interpolation problem for a primary modelling
of fingerprints is to find the \textquotedblleft simplest\textquotedblright\
class of parameter-dependent planar vector fields which cover the widest
possible set of singular points. For the approximating phase portraits given
in Section 4, such a simple class is mainly supposed to cover the Focus, the
Center, the (Improper) Node and the Cusp. An interesting result in this
regards was shown by Andronov et al. \cite{Andronov} and is reported here
from Perko \cite{Perko}:

\bigskip

Let assume that the origin is an isolated singular point of the planar system%
\begin{equation}
\left\{ 
\begin{array}{ccc}
\overset{\cdot }{x} & = & P(x,y) \\ 
\overset{\cdot }{y} & = & Q(x,y)%
\end{array}%
\right.  \label{5}
\end{equation}%
where $P$ and $Q$ are analytic in some neighborhood of the origin. To cover
the Cusp, we consider the case when the Jacobian matrix $J$ at the origin
has two zero eigenvalues, but $J\neq 0.$ Following Andronov et al., system
(5) can be put in the simple form, called normal form,\footnote{%
In fact, at a purely formal level, simpler normal forms can be reached for
nilpotent systems, i.e. systems with a nilpotent linear part (see for
instance Stolovitch \cite{Stolovitch} for a Carleman-linearisation-based
approach). As for the explicit computation of a normal form, it would be
horrible to try to conduct by hand high-order expansions. In this regards, a
Lie-series-based Maple package for symbolic computation of Poincar\'{e}%
-Dulac normal forms in the general case is available at Elsevier's CPC
Library \cite{Zinoun}.}%
\begin{equation}
\left\{ 
\begin{array}{lll}
\overset{\cdot }{x} & = & y \\ 
\overset{\cdot }{y} & = & a_{k}x^{k}(1+h(x))+b_{n}x^{n}y(1+g(x))+y^{2}R(x,y)%
\end{array}%
\right.  \label{6}
\end{equation}%
where $h(x),\ g(x)$ and $R(x,y)$ are analytic in a neighborhood of the
origin, $h(0)=g(0)=0,$ $k\geq 2,$ $a_{k}\neq 0$ and $n\geq 1.$ Then we have
the following:

\begin{theorem}
\cite{Andronov} \ \ Let $k=2m+1$ with $m\geq 1$ in (6) and let $\lambda
=b_{n}^{2}+4(m+1)a_{k}.$ Then if $a_{k}>0$, the origin is a (topological)
Saddle. If $a_{k}<0,$ the origin is (1) a Focus or a Center if $b_{n}=0$ and
also if $b_{n}\neq 0$ and $n>m$ or if $n=m$ and $\lambda <0,$ (2) a Node if $%
b_{n}\neq 0,$ $n$ is an even number and $n<m$ and also if $b_{n}\neq 0,$ $n$
is an even number, $n=m$ and $\lambda \geq 0$ and (3) a critical point with
an elliptic domain if $b_{n}\neq 0,$ $n$ is an odd number and $n<m$ and also
if $b_{n}\neq 0,$ $n$ is an odd number, $n=m$ and $\lambda \geq 0.$\newline
Let $k=2m$ with $m\geq 1$ in (6). Then the origin is (1) a Cusp if $b_{n}=0$
and also if $b_{n}\neq 0$ and $n\geq m$ and (2) a Saddle-Node if $b_{n}\neq
0 $ and $n<m.$
\end{theorem}

\subsubsection{Some normal forms within the ALW classification}

As an interesting case, to compute a normal form for the Elliptical/Circular
Whorl, that is a simple planar system whose phase portrait is a smooth
deformation of the connexion given in Fig.\ref{circular}, we put a Center at
the origin of the $(x,y)$ plane and two face-to-face Cusps at the opposite
points $\pm (1,0);$ then, we compute a vector field of the form (6), for
which the origin and the opposite points are the only singularities, and
whose the respective Taylor polynomials at these points define a local
differential system of the form%
\begin{equation*}
\left\{ 
\begin{array}{lll}
\overset{\cdot }{x} & = & y \\ 
\overset{\cdot }{y} & = & \alpha x+o(x^{2})%
\end{array}%
\right.
\end{equation*}%
for the Center, with $\alpha <0,$ and%
\begin{equation*}
\left\{ 
\begin{array}{lll}
\overset{\cdot }{x} & = & y \\ 
\overset{\cdot }{y} & = & \pm \beta (x\pm 1)^{2}+o((x\pm 1)^{3})%
\end{array}%
\right.
\end{equation*}%
for the Cusps, according to Theorem 6, with $k=2,$ $b_{n}=0,$ $h(x)=0$ and $%
\beta >0.$

\bigskip

\begin{figure}
\begin{center}
\includegraphics{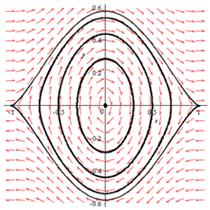} 
\end{center}
\caption{\it Normal form of the
Elliptical/Circular Whorl.}
\label{circular-maple}
\end{figure}

To preserve the center nature of the origin, we only considered the case
where $\overset{\cdot }{y}$ is an univariate function $f(x),$ thus obtaining
symmetric systems with respect to the $x$-axis. Letting $f$ vanish in $0$
with $f^{\prime }(0)=-1$ and in $\pm 1$ with $f^{\prime }(\pm 1)=0,$ $%
f^{\prime \prime }(1)<0$ and $f^{\prime \prime }(-1)>0,$ leads to an
univariate Hermite interpolation problem, for which a solution is given by%
\begin{equation*}
f(x)=-x(x^{2}-1)^{2}
\end{equation*}%
We obtain the normal form 
\begin{equation*}
\left\{ 
\begin{array}{lll}
\overset{\cdot }{x} & = & y \\ 
\overset{\cdot }{y} & = & -x(x^{2}-1)^{2}%
\end{array}%
\right.
\end{equation*}%
whose a Maple-drawn phase portrait is given in Fig.\ref{circular-maple}.

\bigskip

The same reasoning holds for the Spiral where a Focus is to be placed at the
origin instead of a Center. The phase portraits to be connected are
associated to the flows generated by systems of the following form: 
\begin{equation*}
\left\{ 
\begin{array}{lll}
\overset{\cdot }{x} & = & y \\ 
\overset{\cdot }{y} & = & \alpha x+\beta y+o(x^{6})%
\end{array}%
\right.
\end{equation*}%
for the Focus, with $4\alpha +\beta ^{2}<0,$ and%
\begin{equation*}
\left\{ 
\begin{array}{lll}
\overset{\cdot }{x} & = & y \\ 
\overset{\cdot }{y} & = & \pm \gamma (x\pm 1)^{2}+o((x\pm 1)^{3})%
\end{array}%
\right.
\end{equation*}%
for the Cusps, according to Theorem 6, with $k=2,$ $b_{n}=0,$ $h(x)=0$ and $%
\gamma >0.$%
\begin{equation*}
\left\{ 
\begin{array}{lll}
\overset{\cdot }{x} & = & y \\ 
\overset{\cdot }{y} & = & (\gamma +\delta y)(x-1)^{2}+o((x-1)^{3})%
\end{array}%
\right.
\end{equation*}

Fixing $y$ and solving the corresponding univariate Hermite interpolation
problem, with flexible parameters, leads to the normal form%
\begin{equation*}
\left\{ 
\begin{array}{lll}
\overset{\cdot }{x} & = & y \\ 
\overset{\cdot }{y} & = & (y-x/2)(x^{2}-1)^{2}%
\end{array}%
\right.
\end{equation*}%
whose (reversed) phase portrait looks the same as the damped pendulum's in
Fig.\ref{basic}, but with Cusps instead of Saddles.

\bigskip

I could not find, however, a normal form for the configuration in which the
separatrix cycle is approached by the inner spiral as a limit set . In fact,
building a separatrix cycle falls within the scope of the\textit{\ }global
theory of dynamical systems and cannot be dealt with as an interpolation
problem. I proceeded by trial and error within Theorem 6, considering the
Focus as a degenerate singular point $(k=5,$ $b_{n}\neq 0,$ $n=5>m=2,$ $%
h(x)=g(x)=0$ and $\alpha <0)$ of the system%
\begin{equation*}
\left\{ 
\begin{array}{lll}
\overset{\cdot }{x} & = & y \\ 
\overset{\cdot }{y} & = & \alpha (1+y)x^{5}+o(x^{6})%
\end{array}%
\right.
\end{equation*}%
This led to the normal form%
\begin{equation*}
\left\{ 
\begin{array}{lll}
\overset{\cdot }{x} & = & y \\ 
\overset{\cdot }{y} & = & -x^{5}(x^{2}-1)^{2}(1+y(1+x)^{3})%
\end{array}%
\right.
\end{equation*}%
but with no separatrix cycle. Global bifurcations, which are more difficult
to understand than local ones, can be considered within the well-known
Melnikov theory for perturbed planar analytic systems, where parameters
values can be found to characterise bifurcations experienced at homoclinic
or heteroclinic loops. In a more general context, however, we know that
Coppel's problem of determining all possible phase portraits for just a
quadratic planar system and classifying them by means of algebraic
inequalities on the coefficients is insoluble. In other terms, if finding
geometrically all possible connexions - up to a smooth deformation - is
already a difficult task, associating algebraically a normal form to each
connexion is of another kind of difficulty.

\bigskip\ 

For the Twist, where an Improper Node is distributed by two Cusps, and with
the same reasoning as above, the following normal form has been found:%
\begin{equation*}
\left\{ 
\begin{array}{lll}
\overset{\cdot }{x} & = & y \\ 
\overset{\cdot }{y} & = & (2y-x)(x-1)^{2}%
\end{array}%
\right.
\end{equation*}%
and for Plain or Tented Arches, the straight or tented curves can be merely
seen as a smooth deformation of the phase portrait of the flow generated by
a trivial system with zero $\overset{\cdot }{y}$ and constant $\overset{%
\cdot }{x},$ or, as the case may be, a simple system corresponding to a
topsy-turvy Cusp, as shown in Fig.\ref{arch}. Finally, the lack of normal
forms for Loops and for some special cases of Whorls should be noted, a
question to be sought not only from a traditional planar interpolation
viewpoint, but also involving infinite interpolation points or, eventually,
considering the qualitative behaviour on invariant planes of high
dimensional dynamical systems.\footnote{%
The increment in dimension, however, doesn't necessarily mean increment in
behavioural complexity. For example, no strange attractor is required in
fingerprints modelling, and until proved otherwise, no individual has been
identified with the Lorenz butterfly on the thumb!}

\subsection{Patched phase portraits, matched fingerprints:\newline
Some perspective}

\begin{figure}
\begin{center}
\includegraphics{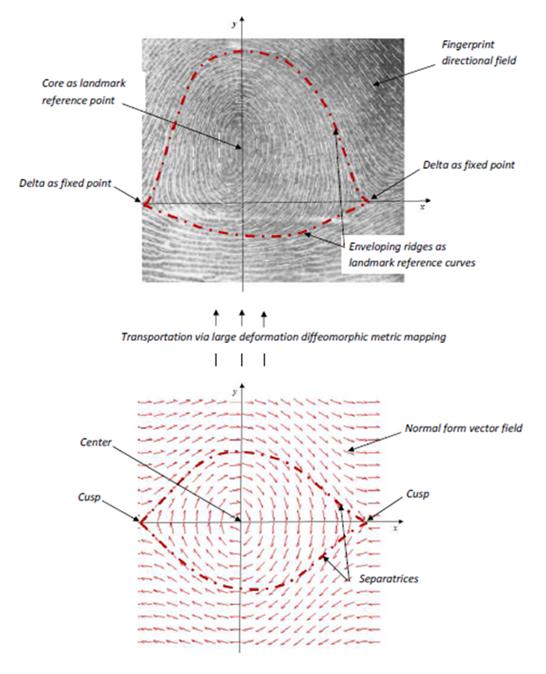} 
\end{center}
\caption{\it Landmark matching
problem for the Elliptical/Circular Whorl.}
\label{LDDMM}
\end{figure}

Now that local phase portraits have been patched together within a
structure-preserving normal form connexion, that is a connexion preserving
the nature of singular points, but not necessarily their position, as for
example in the case of the Elliptical/Circular Whorl, the problem is to find
the simplest deformation allowing transportation from the normalised phase
portrait as a source image to the underlying fingerprint's stream of ridges
as a target image (Fig.\ref{LDDMM}). This idea, from the academic discipline
of computational anatomy, will be developped in future work. Specifically,
instead of roughly considering smooth deformations, we will be using the
terminology of large deformation diffeomorphic metric mapping\textit{\ }%
(LDDMM), as follows from earlier works by Christensen et al. \cite{Chris}
and Trouv\'{e} \cite{Trouve} on pattern matching in image analysis. An
intrinsic Cartesian coordinate system will be then associated to each
category of fingerprints, where, in the case of a whorled pattern for
example, Deltas are maintained fixed along the $x$-axis and the core, say a
Center, is to be displaced from the origin to its real position along the $y$%
-axis via a LDDMM. Besides Deltas and cores, other reference points (and,
more generally, portions of curves) can be considered within the landmark
matching problem, mainly those corresponding to the maximal curvature of the
enveloping ridges of the pattern, i.e. mathematically speaking, the common
separatrices of the Cusps. Many LDDMM algorithmic codes are available in the
literature and should prove valuable for such a purpose, but whichever
software suite is being used, the lesser the input data the better is
modelling, as diffeomorphic mapping parameters have to be added to those of
the normal form to optimally encode the fingerprint.

\bigskip

Even more to be welcomed, though difficult to express in terms of elementary
functions, then to implement, is to consider nonautonomous differential
systems, normal form source image and fingerprint's orientation target image
being seen as screenshots of a moving nonlinear phase portrait. To encode
the targeted orientation field, only model's coefficients and exact time of
simulation are required. It is then no violent misuse of metaphor to
compare, for example, a moving Focus of a nonautonomous dynamical system to
a moving cyclone in a meteorological forecast map. Interesting, but easier
said than done!

\section{Conclusion: \textit{Elementary, my dear Galton}}

\begin{figure}
\begin{center}
\includegraphics{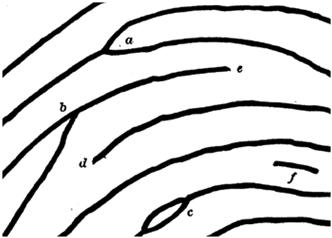} 
\end{center}
\caption{\it Characteristic
peculiarities in ridges, about 8 times the natural size, as captured from
Plate 3 of Galton's book.}
\label{detail}
\end{figure}


Let us go back now to the main question stated at header level: can a
fingerprint be modelled by a differential equation? Obviously, if that is
meant to be a model which faithfully reproduces the fingerprint, right down
to the finest detail about the papillary ridges (Fig.\ref{detail}), the
answer is certainly a no. In fact, according to Definition 4.1, such a
modelling should ideally be starting from a finite set of initial
conditions, i.e. a set of reference points in the pattern, as for example
those where ridges stop abruptly, then, by some fundamental law,\footnote{%
This amounts to writing down one's unique genetic code as a differential
equation!} reproducing exactly the stream of ridges of the fingerprint
under study. We already know, however, that if solving differential
equations in terms of elementary functions is generally impossible,
conversely, trying to assign an explicitly written differential equation to
an observed system of curves is most often illusory. Put another way, the
family album of our elementary functions is too small to allow closed-form
expressions to be derived for any curve encountered in nature.\footnote{%
Like smooth functions, which in terms of categories have been shown by
Banach to be negligible for the class of continuous functions, i.e. almost
all continuous functions are everywhere nondifferentiable, curves that can
be expressed in terms of elementary functions, allowing special ones and any
combination of them, should prove exotic among a collection of curves taken
for example from a child's drawing!} As an approximation, following the
work by Ford \cite{Ford} on visualisation in fluid dynamics, the least
squares method can be used to adjust the coefficients of a polynomial model
vector fields to minimise the difference between its integral curves and the
observed flow streamline. Although practically interesting, at a purely
mathematical level, some hidden problems from dynamical systems theory
cannot however come to light if the modelling is entrusted to a mere
least-squares-based approach. For example, as a first step for a primary
classification, the choice of a diffeomorphic conjugation of the fingerprint
directional field as a symmetric normalised vector field was guided by an
intention to keep under control and mitigate the risk (of changing in
nature) faced by nonhyperbolic singular points when proceeding to a
connexion whose linearisation escapes the control of the Hartman-Grobman
theorem. Of course, a symmetry condition is too restrictive for a vector
field and is poorly adapted to the nonsymmetric character trait of a
fingerprint orientation field. However, on the one hand, and as said by the
author in \cite{Poincare2}, necessary and sufficient Poincar\'{e}'s
conditions to ensure the preservation of the center-like nature of a
singular point are too cumbersome so that it would be difficult to proceed
to modelling within an asymmetric framework. On the other hand, it should
not be overlooked that, basically, the symmetric normal form is just an
example of a starting point in the modelling process, a stepping-stone which
seems in every respect preferable, but any structure-preserving connexion
could be safely used as a source image to be diffeomorphically carried
towards the targeted fingerprint orientation image.

\bigskip

Another framework within which the present work could have been conducted
would be a reading where curves under study are considered as contour lines
of a Poincar\'{e}'s topographic system,\footnote{%
As already pointed out by Galton in [1, p.77], \textit{a curious optical
effect is connected with the circular forms, which becomes almost annoying
when many specimens are examined in succession. They seem to be cones
standing bodily out from the paper. This singular appearance becomes still
more marked when they are viewed with only one eye; no stereoscopic guidance
then correcting the illusion of their being contour lines.}} double-points
of compound cycles being assigned to terrain's passes and singular points
not belonging to a cycle corresponding to bottoms and summits. This
however should lead to the same order of difficulty as discussed above, and
in fact, building an algebraic relief surface whose projected contour lines
form a carbon copy of the fingerprint pattern would be as complicated as
writing down a differential system whose solution curves align with the
fingerprint stream of ridges.

\bigskip

Finally, beyond biometry, forensic science, or the question of whether the
qualitative theory of differential equations should be taught to budding
Sherlock Holmes, the main beneficiary of the present study is perhaps
mathematics itself. In fact, as the reader may have guessed, fingerprints
(and any oriented texture in general) were merely a pretext for raising
interesting questions in dynamical systems theory, as for example, the
rather formidable connexion problem for a set of randomly distributed and
natured singular points in space.

\newpage

\begin{figure}
\begin{center}
\includegraphics{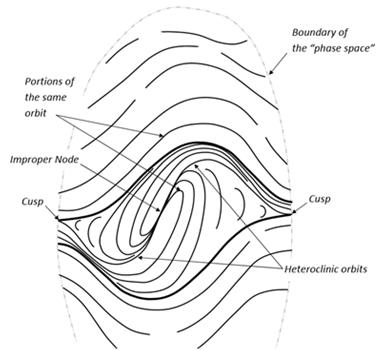} 
\end{center}
\end{figure}


\begin{center}
\bigskip
\end{center}

\textit{When I look at the Twist, I almost feel like singing Chubby Checker's%
}

\textit{\textquotedblleft Let's Twist Again\textquotedblright\ from the
'60s, with some modifications in the lyrics .. }

\bigskip

\begin{center}
\textit{Come on everybody !}

\textit{Print your hands !}

\textit{(Inky) All you looking good !}

\textit{I'm gonna (rock'n) roll your thumb}

\textit{It won't take long !}

\textit{We're gonna print the Twist}

\textit{And it goes like this :}

\bigskip

\textit{Come on let's twist again,}

\textit{Like it is on the last phalanx !}

\textit{Yeaaah, let's twist again,}

\textit{Like we did for the last Whorl !}

\bigskip

(...)

\newpage
\end{center}

\end{document}